\documentclass[11pt]{amsart}

\evensidemargin\oddsidemargin
\begin{document}

\baselineskip=18pt
\setcounter{page}{1}
    
\newtheorem{Thm}{Th\'eor\`eme 1\!\!}
\newtheorem{Rqs}{Remarques 2\!\!}
\newtheorem{Thmm}{Th\'eor\`eme 3\!\!}
\newtheorem{Cor}{Corollaire 4\!\!}
\newtheorem{Corr}{Corollaire 5\!\!}
\newtheorem{Prp}{Proposition 6\!\!}
\newtheorem{Rq}{Remarque 7\!\!}

\renewcommand{\theRqs}{}
\renewcommand{\theThm}{}
\renewcommand{\theThmm}{}
\renewcommand{\theCor}{}
\renewcommand{\theCorr}{}
\renewcommand{\thePrp}{}
\renewcommand{\theRq}{}

\def\a{\alpha}
\def\CC{{\mathbb{C}}} 
\def\Ea{E_\a}
\def\EE{{\mathbb{E}}} 
\def\elaw{\stackrel{d}{=}}
\def\eps{\varepsilon}
\def\hS{{\hat S}}
\def\hT{{\hat T}}
\def\hX{{\hat X}}
\def\ii{{\rm i}}
\def\lb{\lambda}
\def\lacc{\left\{}
\def\lcr{\left[}
\def\lpa{\left(}
\def\lva{\left|}
\def\NN{{\mathbb{N}}} 
\def\pb{{\mathbb{P}}}
\def\rl{{\mathbb{R}}}
\def\racc{\right\}}
\def\rcr{\right]}
\def\rpa{\right)}
\def\Un{{\bf 1}}
       
\def\abstractname{R\'esum\'e}

\def\keywordsname{Mots-cl\'es}

\def\subjclassname{Classification AMS}

\bibliographystyle{french}

\newcommand{\fin}{\vspace{-0.6cm}
                  \begin{flushright}
                  \mbox{$\Box$}
                  \end{flushright}
                  \noindent}

\title[Fonctions de Mittag-Leffler et processus stables]{Fonctions de Mittag-Leffler et processus de L\'evy stables sans saut n\'egatif}

\author[Thomas Simon]{Thomas Simon}

\address{Laboratoire Paul Painlev\'e, U. F. R. de Math\'ematiques, Universit\'e de Lille 1, F-59655 Villeneuve d'Ascq Cedex. {\em Adresse \'electronique} : {\tt simon@math.univ-lille1.fr}}



\begin{abstract} On remarque que la fonction $\Ea(x^\a) - \a x^{\a -1} \Ea'(x^\a)$ est compl\`etement monotone pour tout $\a\in [1,2].$ Gr\^ace \`a l'expression de sa densit\'e de Bernstein, on en d\'eduit une identit\'e en loi entre suprema de processus de L\'evy $\a-$stables compl\`etement asym\'etriques. Dans le cas spectralement positif, on retrouve l'expression d'une constante de petites d\'eviations unilat\`eres qui avait \'et\'e obtenue avec une autre m\'ethode par Bernyk, Dalang et Peskir \cite{BDP}.
\end{abstract}

\maketitle

\section{Une propri\'et\'e des fonctions de Mittag-Leffler}

La fonction de Mittag Leffler $\Ea$ d'indice $\a > 0$ est d\'efinie pour tout $z\in \CC$ par le d\'eveloppement en s\'erie enti\`ere
\begin{equation}
\label{Mit1}
\Ea(z)\; =\; \sum_{n=0}^{\infty} \frac{z^n}{\Gamma (1 +\a n)}
\end{equation}
o\`u $\Gamma$ d\'esigne la fonction Gamma d'Euler. En \'ecrivant $1/\Gamma$
comme une int\'egrale curviligne le  long d'un chemin de Hankel - voir par
exemple \cite{D} Chapitre IX.4 et Figure 56 :
$$\frac{1}{\Gamma (1 +\a n)}\; =\; \frac{1}{2\pi \ii}\int_H t^{-\a n -1} e^t \, dt,$$
et en int\'egrant la s\'erie terme \`a terme on trouve la repr\'esentation bien connue
\begin{equation}
\label{Mit2}
\Ea(z)\; =\; \frac{1}{2\pi \ii}\int_{H_z}\frac{t^{\a -1}e^t}{t^\a - z} dt
\end{equation}
o\`u $H_z$ d\'esigne un chemin de Hankel contournant le cercle centr\'e en l'origine de rayon $\vert z\vert.$ En particulier on a pour tout $x\ge 0$
$$\Ea(x^\a)\; =\; \frac{1}{2\pi \ii}\int_{H_{x^\a}} \frac{t^{\a -1}e^t }{t^\a - x^\a}dt.$$
On peut appliquer la m\^eme proc\'edure \`a la fonction suivante d\'efinie sur $(0, +\infty):$
$$x\mapsto  \a x^{\a -1} \Ea'(x^\a)\; =\; \sum_{n=1}^{\infty} \frac{x^{\a n -1}}{\Gamma (\a n)},$$
pour trouver la repr\'esentation
$$\a x^{\a -1} \Ea'(x^\a)\; =\; \frac{1}{2\pi \ii}\int_{H_{x^\a}} \frac{x^{\a -1}e^t }{t^\a - x^\a} dt.$$
On d\'efinit maintenant la fonction 
$$D_\a (x) \; =\; \Ea(x^\a) - \a x^{\a -1} \Ea'(x^\a)\; =\; \frac{1}{2\pi \ii}\int_{H_{x^\a}} e^t \lpa\frac{t^{\a -1} - x^{\a -1}}{t^\a - x^\a}\rpa dt$$
pour tout $x > 0$ et on la prolonge par continuit\'e en z\'ero en posant $D_\a
(0) = 1.$ On rappelle qu'une fonction $f : (0, +\infty)\to \rl$ est
compl\`etement monotone - ce que nous noterons par CM - si elle est ind\'efiniment d\'erivable et si
$$(-1)^n f^{(n)}(x) \; \ge \; 0$$
pour tout $n\in\NN, x > 0.$ Le th\'eor\`eme de Bernstein - voir par exemple
\cite{F} Chapitre XIII.4 - dit que si $f$ est CM et si $f(x) \to 1$ quand $x\to 0,$ alors il existe une unique mesure de probabilit\'e $\mu$ sur $\rl^+$ telle que
$$f(x) \; =\; \int_0^\infty e^{-xt} \mu(dt)$$
pour tout $x > 0,$ la r\'eciproque \'etant imm\'ediate. 

\begin{Thm} Pour tout $\alpha\in [1,2],$ la fonction $D_\a$ est CM.
\end{Thm}

\noindent
{\em Preuve} : La propri\'et\'e est \'evidente pour $\alpha = 1$ puisque $D_1(x) = 0$ et pour $\alpha = 2$ puisque l'on voit facilement que $D_2(x) = e^{-x},$ prototype d'une fonction CM avec mesure de Bernstein $\mu = \delta_1.$ Supposons maintenant $\alpha \in (1,2)$ et fixons $x > 0.$ En raisonnant comme dans \cite{D} Chapitre VII.7, la fonction
$$g : t \mapsto \frac{t^{\a -1} - x^{\a -1}}{t^\a - x^\a}$$
se prolonge analytiquement dans tout le plan fendu $\CC/(-\infty, 0]$ en ayant pos\'e $g(x^\a) = (1 - 1/\a)/x$. Par le th\'eor\`eme de Cauchy, on a donc
$$D_\a (x) \; =\; \frac{1}{2\pi \ii}\int_H e^t \lpa\frac{t^{\a -1} - x^{\a -1}}{t^\a - x^\a}\rpa dt$$
o\`u $H$ est un chemin de Hankel entourant l'origine ind\'ependamment de $x$. Pour tout $\rho > 0,$ on d\'ecompose classiquement $H$ en la branche inf\'erieure de $\rl^-$ parcourue de $-\infty$ \`a $-\rho$, le cercle centr\'e en l'origine de rayon $\rho$ parcouru dans le sens trigonom\'etrique, et la branche sup\'erieure de $\rl^-$ parcourue de $-\rho$ \`a $-\infty.$ Comme $x >0,$ l'int\'egrale le long du cercle tend vers 0 quand $\rho \to 0.$ L'int\'egrale le long des deux branches vaut d'autre part

\begin{eqnarray*}
I_\rho & = & \frac{1}{2\pi \ii}\lpa  \int_\rho^\infty e^{-s} \lpa \frac{e^{-\ii\pi (\a -1)}s^{\a -1} - x^{\a -1}}{e^{-\ii \pi\a}s^\a  - x^\a}\rpa ds \; -\; \int_\rho^\infty e^{-s} \lpa \frac{e^{\ii \pi (\a -1)}s^{\a -1} - x^{\a -1}}{e^{\ii \pi\a}s^\a  - x^\a}\rpa ds\rpa
\\
& \to & \frac{1}{2\pi \ii}\lpa  \int_0^\infty e^{-s} 
\lpa \frac{e^{\ii \pi \a }s^{\a -1} + x^{\a -1}}{e^{\ii \pi\a}s^\a  - x^\a} \; -\; \frac{e^{-\ii \pi \a}s^{\a -1} + x^{\a -1}}{e^{-\ii \pi\a}s^\a  - x^\a}\rpa ds\rpa
\end{eqnarray*}
quand $\rho \to 0.$ Apr\`es changement de variable $s = xu,$ on trouve

\begin{eqnarray*}
D_\a (x) & = & \frac{1}{2\pi \ii} \int_0^\infty e^{-xu} 
\lpa \frac{e^{\ii \pi \a }u^{\a -1} + 1}{e^{\ii \pi\a}u^\a  - 1}\; -\; \frac{e^{-\ii \pi \a}u^{\a -1} + 1}{e^{-\ii \pi\a}u^\a  - 1}\rpa du \\
& = & \frac{-\sin \pi\a}{\pi}\int_0^\infty e^{-xu} \lpa \frac{u^{\a
    -1}(1+u)}{u^{2\a} - 2 u^\a\cos \pi\a + 1}\rpa du,
\end{eqnarray*}
ce qui entra\^{\i}ne que $D_\a$ est CM avec pour mesure de Bernstein
\begin{equation}
\label{Brn}
\mu_\a(dt) \; = \; \frac{(-\sin \pi\a)t^{\a -1}(1+t)}{\pi(t^{2\a} - 2t^\a\cos \pi\a +1)} dt.
\end{equation}
\fin

\begin{Rqs}{\em (a) La question de la monotonicit\'e compl\`ete ne se pose
    \'evidemment pas pour $\Ea$ elle-m\^eme mais il est classique \cite{P1}
    que $x\mapsto \Ea(-x)$ est CM pour tout $\a\in
    (0,1].$ Cette propri\'et\'e avait \'et\'e obtenue auparavant par Feller avec un
    argument probabiliste et nous y reviendrons dans la section suivante. Nous
    renvoyons \`a \cite{Sc} pour des r\'esultats plus g\'en\'eraux, qui
    entra\^{\i}nent en particulier que $x\mapsto \Ea(-x)$ n'est pas CM quand $\a > 1.$

\vspace{2mm}

\noindent
(b) On voit facilement que $\mu_\a \Rightarrow 0$ quand $\alpha \to 1$ et on peut montrer que 
$\mu_\a \Rightarrow \delta_1$ quand $\alpha \to 2.$ Jean-Fran\c{c}ois Burnol
m'a aussi fait remarquer que la fonction $\a\mapsto \Ea(x)$ est en fait analytique
dans un c\^one contenant $\rl^+$ pour tout $x > 0.$ Quand $\a\in ]0,1[$ les calculs
pr\'ec\'edents restent valables et entra\^{\i}nent que $-D_\a$ est
CM avec mesure de Bernstein
$$\frac{(\sin \pi\a)t^{\a -1}(1+t)}{\pi(t^{2\a} - 2t^\a\cos \pi\a +1)} dt.$$
En revanche, quand $\a > 2$ il semble que l'on perde
toute monotonicit\'e puisque la fonction $g$ ci-dessus admet alors au moins deux
p\^oles \`a l'int\'erieur du contour. On peut par exemple calculer
$$D_4(x) \; =\; \frac{1}{2}(e^{-x} + \cos x + \sin x).$$
Il n'est pas impossible qu'il faille retrancher ou ajouter des d\'eriv\'ees
successives \`a $D_\a$ pour retomber sur une fonction CM.}

\end{Rqs}

Consid\'erons maintenant pour tout $\a\in [1,2]$ la fonction $F_\a$ d\'efinie sur $\rl^+$ par
$$F_\a(x) \; = \; D_\a(x^{1/\a})\; =\; \Ea(x) - \a x^{1-1/\a}\Ea'(x).$$
Comme $x\mapsto x^{1/\a}$ est une fonction positive de 
d\'eriv\'ee CM, par le Th\'eor\`eme 1 on sait en appliquant le Crit\`ere 2 dans \cite{F} Chapitre XIII. 4 que $F_\a$ est \'egalement CM. En fait on aurait pu obtenir ce r\'esultat, qui est l\'eg\`erement plus faible que le Th\'eor\`eme 1, par un argument probabiliste. C'est ce que nous allons commencer par expliquer dans la section suivante.

\section{Une identit\'e en loi pour les processus stables}

On fixe $\a\in (1, 2]$ et on consid\`ere un processus de L\'evy $\a-$stable
sans saut n\'egatif $X = \lacc X_t, \; t\ge 0\racc,$ normalis\'e tel que
$\EE[e^{-X_1}] = e.$ On pose $\hX = -X$ qui est un processus de L\'evy
$\a-$stable sans saut positif et l'on renvoie aux chapitres VII et VIII de \cite{B} pour plus de d\'etails concernant $\hX$. La formule de L\'evy-Khintchine s'\'ecrit
$$\EE[e^{-\lb X_t}]\; =\; \EE[e^{\lb \hX_t}]\; =\; e^{t\lb^a}
$$ 
pour tout $t, \lb \ge 0,$ de sorte qu'avec cette normalisation, quand $\alpha = 2$ on a $X\elaw\hX\elaw \sqrt{2} B$ o\`u $B$ est un mouvement brownien standard. On d\'efinit les suprema
$$S_t \; =\; \sup\{X_s, \; s\le t\}\quad \mbox{et}\quad\hS_t \; =\; \sup\{\hX_s, \; s\le t\},\quad t\ge 0.$$ Pour tout $q > 0$, soit $\tau_q$ une variable exponentielle de param\`etre $q$ ind\'ependante de $X$. Comme $S_t = - \inf\{\hX_s, \; s\le t\}$ pour tout $t$, la factorisation de Wiener-Hopf pour $\hX$ - voir \cite{B} Chapitre VII.1 Formule (3) - entra\^{\i}ne que pour tout $\lb > q^{1/\a} > 0$
\begin{equation}
\label{WH}
\EE[e^{-\lb S_{\tau_q}}]\; =\; \frac{q(\lb - q^{1/\a} )}{q^{1/\a}(\lb^\a - q)}
\end{equation}
et apr\`es int\'egration par parties, il vient
$$\int_0^\infty e^{-\lb x}\pb[S_{\tau_q} \ge x]\, dx\; =\; \frac{\lb^{\a -1}}{\lb^\a - q}  - \frac{q^{1-1/\a}}{\lb^\a - q} $$
pour tout $\lb > q^{1/\a} > 0$. En fait les deux formules ci-dessus sont valables pour tout $\lb, q >0$ avec un prolongement par continuit\'e imm\'ediat en $\lb = q^{1/\a}.$

Le lien avec les fonctions de Mittag-Leffler est le suivant : en int\'egrant (\ref{Mit1}) terme \`a terme on peut aussi obtenir l'expression de la transform\'ee de Laplace de $\Ea(q x^\a)$ :
 pour tout $\lb > q^{1/\a}$ on a 
\begin{equation}
\label{Mit3}
\int_0^\infty e^{-\lb x} \Ea( qx^\a)\, dx \; =\; \frac{\lb^{\a -1}}{\lb^\a - q}
\end{equation}
et en int\'egrant par parties, il vient
\begin{eqnarray*}
\int_0^\infty e^{-\lb x} (\Ea( qx^\a) - \a (qx^\a)^{1-1/\a} \Ea'( qx^\a))\, dx  & = & \frac{\lb^{\a -1}}{\lb^\a - q} - \frac{q^{1-1/\a}}{\lb^\a - q}\\
& = & \int_0^\infty e^{-\lb x}\pb[S_{\tau_q} \ge x]\, dx.
\end{eqnarray*}
Signalons que la formule (\ref{Mit3}) avait d\'ej\`a \'et\'e utilis\'ee dans
\cite{H} pour obtenir de nouvelles correspondances symboliques entre fonctions
$\Ea$ d'indices diff\'erents, et plus r\'ecemment dans \cite{B1} pour calculer une constante de petites d\'eviations en norme uniforme pour $X$. L'inversion de la transform\'ee de Laplace donne 
\begin{equation}
\label{L1}
\pb[S_{\tau_q} \ge x]\; = \; \Ea( qx^\a) - \a (qx^\a)^{1 -1/\a} \Ea'( qx^\a)
\end{equation}
pour tout $q, x >0.$ Il n'est pas \'evident a priori que la fonction $x
\mapsto \pb[S_{\tau_1} \ge x]$ soit compl\`etement monotone. En revanche, on
peut retrouver \`a l'aide de (\ref{L1}) le r\'esultat \'enonc\'e \`a la fin de
la premi\`ere section. Pour cela, on introduit pour tout $x > 0$ le premier temps de sortie
$$T_x = \inf\{ t > 0, \; X_t > x\}.$$ 
Comme $X$ est auto-similaire d'indice $1/\a$ on a $T_x \elaw x^\a T_1,$ et la continuit\'e \`a droite de $X$ entra\^{\i}ne d'autre part $\{S_t \ge x\} = \{ T_x \le t\}$ p.s. Apr\`es changement de variable $y = qx^\a$ dans (\ref{L1}), on trouve
\begin{equation}
\label{L2}
\int_0^\infty e^{-yt} \pb[T_1 \in dt]\; =\; \Ea(y) - \a y^{1-1/\a} \Ea'(y), \qquad y > 0,
\end{equation}
et ceci entra\^{\i}ne que la fonction \`a droite est compl\`etement monotone avec mesure de Bernstein $\pb[T_1\in dt].$ On va maintenant inverser (\ref{L2}) et donner une expression de la loi de $T_1$. Pour $\a = 2$ et $\hX \elaw \sqrt{2} B,$ le terme de droite vaut $e^{-\sqrt{x}}$ et on retrouve l'expression classique pour la densit\'e de $T_1$ qui est $e^{-1/4t}/2t\sqrt{\pi t}.$ On fixe maintenant $\a\in ]1,2[$ et on note 
$$\hT_1 \; = \; \inf\{ t > 0, \; \hX_t > 1\}.$$
Il est bien connu et facile \`a voir - voir le Th\'eor\`eme VII.1 dans \cite{B} - que $\EE[e^{-\lb\hT_1}] = e^{-\lb^{1/\a}}$ pour tout $\lb > 0$, autrement dit que $\hT_1$ est une variable stable positive d'indice $1/\a.$ Sa densit\'e $f_{\hT_1}$ peut \^etre donn\'ee sous forme int\'egrale - voir \cite{P2} pour le calcul :
$$f_{\hT_1}(t)\; =\; \frac{1}{\pi}\int_0^\infty e^{-(tu+u^{1/\a}\cos(\pi/\a))}\sin(u^{1/\a}\sin(\pi/\a))\, du,$$
expression qui se simplifie l\'eg\`erement pour $\a = 3/2$ en termes de la fonction de Whittaker $W_{-1/2, -1/6}$ - voir \`a nouveau \cite{P2}. On consid\`ere d'autre part une variable $T$ ind\'ependante de $\hT_1$ et de densit\'e
$$f_T(t)\; =\;\frac{(-\sin \pi\a)(1+t^{1/\a})}{\pi\a(t^2 - 2t \cos \pi\a  +1)}$$
sur $\rl^+.$ Il est \'evident apr\`es un changement de variable dans
(\ref{Brn}) que $T$ est bien une variable al\'eatoire finie p.s. mais je n'ai
pas trouv\'e d'exemples o\`u $T$ jou\^at un r\^ole dans la litt\'erature. En
particulier, malgr\'e les apparences, il ne semble pas que $T$ soit 
reli\'ee \`a une loi de l'arcsinus - voir \cite{B} Th\'eor\`eme III.6 - ou
encore \`a un noyau fractionnaire de type Kato - voir \cite{Y} Formule XI.11.6, puisque $\a > 1.$

\begin{Thmm} Avec les notations pr\'ec\'edentes, on a pour tout $\a\in ]1, 2[$
$$T_1\;\elaw\; T\,\times\,\hT_1.$$
\end{Thmm}
\noindent 
{\em Preuve} : En combinant (\ref{L2}) et le th\'eor\`eme 1, on a pour tout $x >0$
$$\int_0^\infty e^{-x^\a t} \pb[T_1 \in dt]\; =\; D_\a(x) \; =\; \int_0^\infty e^{-xt} \mu_\a(dt)$$
o\`u l'on a repris la notation (\ref{Brn}). D'autre part, comme par d\'efinition
$$e^{-xt}\; =\; \int_0^\infty e^{-(xt)^\a u} f_{\hT_1}(u)\; du,$$
on d\'eduit par Fubini et en faisant le changement de variable $s = t^\a u$
\begin{eqnarray*}
\int_0^\infty e^{-x^\a t} \pb[T_1 \in dt] & = & \int_0^\infty \int_0^\infty e^{-(xt)^\a u}\lpa  \frac{(-\sin \pi\a) t^{\a -1}(1+t)}{\pi(t^{2\a} - 2t^\a\cos \pi\a  +1)} \rpa f_{\hT_1}(u) du\, dt \\
& = &  \int_0^\infty  e^{-x^\a s}\lpa\int_0^\infty f_{\hT_1}(u) f_T(s/u)\frac{du}{u}\rpa\, ds, \\
\end{eqnarray*}
ce qui prouve en inversant la transform\'ee de Laplace que $T_1$ a pour densit\'e
$$f_{T_1}(t)\; =\; \int_0^\infty f_{\hT_1}(u) f_T(t/u)\frac{du}{u}\; =\; \int_0^\infty f_{T}(u) f_{\hT_1}(t/u)\frac{du}{u}$$
qui est bien celle de $T\times\hT_1.$

\fin 

Cette identit\'e en loi peut se d\'ecliner de plusieurs fa\c{c}ons. On sait par exemple par auto-similarit\'e que $T_1\elaw 
S_1^{-\a}$ et $\hT_1\elaw \hS_1^{-\a},$ d'o\`u
\begin{equation}
\label{hS}
S_1\; \elaw\; T^{-1/\a}\times\, \hS_1
\end{equation}
avec les notations pr\'ec\'edentes. Ceci donne une autre relation entre $S_1$ et les fonctions de Mittag-Leffler que celle donn\'ee par (\ref{L1}). On sait en effet  que $\hS_1$ suit une {\em loi} de Mittag-Leffler d'indice $1/\a$ au sens o\`u
$$\EE[e^{-\lb \hS_1}]\; =\; E_{1/\a} (-\lb), \qquad \lb \ge 
0$$
et l'on renvoie pour cela \`a l'Exercice 29.18 dans \cite{S} qui
donne une preuve probabiliste du r\'esultat de Pollard \cite{P1}. En
utilisant (\ref{hS}) et un changement de variable que nous laissons au
lecteur, on obtient donc une relation probablement inutile, mais in\'edite \`a
notre connaissance, entre fonctions de Mittag-Leffler d'indices diff\'erents :

\begin{Cor} On fixe $\a\in ]1,2[$ et on note $f_{1/\a}$ la densit\'e de Bernstein associ\'ee \`a la fonction compl\`etement monotone $x\mapsto E_{1/\a} (- x)$. Avec les notations pr\'ec\'edentes, on a
$$D_\a(x) \; =\; \int_0^\infty\int_0^\infty f_{1/\a} (u) e^{-sx^\a/u^\a}\frac{(-\sin \pi\a)(1+s^{1/\a})}{\pi\a(s^2 - 2s\cos \pi\a +1)} duds$$ 
pour tout $x > 0.$ 
\end{Cor}
 A l'aide du th\'eor\`eme de Skorohod qui stipule que $\hS_1\elaw\hX_1$ sachant $\hX_1 > 0$ - voir pour cela les Exercices 29.7 et 29.18 dans \cite{S}, on obtient enfin une relation entre les lois de $S_1$ et $X_1$:
 
\begin{Corr} Avec les notations pr\'ec\'edentes on a
$$S_1\; \elaw\; -(T^{-1/\a}\times\, X_1)\;\; \mbox{{\rm sachant}}\;\; X_1 < 0.$$
\end{Corr}
 
\section{Un calcul de constante}

On va s'int\'eresser maintenant au comportement asymptotique de la fonction de r\'epartition de $T_1$ en l'infini. On sait par la Proposition VIII.2 dans \cite{B} et l'identit\'e $T_1\elaw S_1^{-\a}$ que
\begin{equation}
\label{Kpp}
\pb[T\ge t]\; \sim\; \kappa t^{1/\a -1}, \qquad t\to +\infty
\end{equation}
o\`u $\kappa$ est une constante d\'ependant de la loi de la trajectoire du processus $X$ - plus pr\'ecis\'ement de celle de son processus d'\'echelle bivari\'e $(L^{-1}, H)$ en suivant les notations de \cite{B}. Les r\'esultats pr\'ec\'edents permettent de calculer $\kappa$ :

\begin{Prp} La constante $\kappa$ dans {\em (\ref{Kpp})} vaut $1/\Gamma(\a)\Gamma(1/\a).$
\end{Prp}

\noindent
{\em Preuve :} Dans le cas $\alpha = 2$ et $\kappa = 1/\sqrt{\pi},$ le r\'esultat est bien connu et d\'ecoule directement de l'expression de la densit\'e de $T_1$ rappel\'ee plus haut. Dans le cas $\alpha \in ]1,2[,$ on d\'eduit d'abord de (\ref{L2}) et du Th\'eor\`eme 1 que pour tout $\lb > 0$
$$\EE[e^{-\lb^\a T_1}]\; =\; \frac{-\sin \pi\a}{\pi}\int_0^\infty e^{-\lb u} \lpa \frac{u^{\a -1}(1+u)}{u^{2\a} - 2u^\a\cos \pi\a  + 1}\rpa du.$$
Apr\`es une int\'egration par parties et un changement de variables, ceci implique
$$\lb^\a \int_0^\infty e^{-\lb^\a t} \pb [T_1\ge t]\, dt\; =\; 
  \lb \int_0^\infty e^{-\lb t} \pb [T\ge t^\a]\, dt$$
avec les notations pr\'ec\'edentes pour $T$. Gr\^ace \`a l'expression de $f_T$, on sait d'autre part que
$$\pb [T\ge t^\a]\; \sim\; \frac{\sin \pi\a}{\pi(1-\a)}t^{1-\a}, \qquad t\to \infty.$$
En utilisant deux fois de suite un th\'eor\`eme taub\'erien - celui donn\'e par exemple dans \cite{F} Chapitre XIII.5, Th\'eor\`eme 4, on voit que
$$\pb[T\ge t]\; \sim\; \frac{\sin \pi\a \Gamma (2-\a)}{\pi(1-\a)\Gamma(1/\a)} t^{1/\a -1}\; =\; \frac{1}{\Gamma(\a)\Gamma(1/\a)} t^{1/\a -1}, \qquad t\to +\infty,$$
o\`u dans la derni\`ere \'egalit\'e on a utilis\'e successivement la formule de r\'ecurrence et la formule des compl\'ements pour la fonction Gamma.

\fin

\begin{Rq} {\em En \'ecrivant
$$\lb^\a \int_0^\infty e^{-\lb^\a t} \pb [T_1\le t]\, dt\; =\; 
  \lb \int_0^\infty e^{-\lb t} \pb [T\le t^\a]\, dt, \quad \lb > 0,$$
on voit \`a cause de l'asymptotique $\pb [T\le t^\a]\sim(-\sin \pi\a /\pi\a)t^\a$ et \`a nouveau par un th\'eor\`eme taub\'erien que
$$\pb[T_1\le t]\; \sim\; \frac{1}{\pi} (- \Gamma(\a)\sin \pi\a) t\; = \; \lpa\frac{-1}{\Gamma(1-\a)}\rpa t, \qquad t\to 0^+,$$
r\'esultat d\'ej\`a connu en combinant l'identit\'e $T_1 \elaw S_1^{-\a}$ avec
la Proposition VIII.4 de \cite{B} et l'estim\'ee (14.37) de \cite{S} dans le
cas particulier $\rho = 1- 1/\a.$}
\end{Rq} 

\section{Remarques finales}

En fait le r\'esultat de la Proposition 6 avait lui aussi d\'ej\`a \'et\'e obtenu,
plus r\'ecemment et avec une normalisation diff\'erente,  dans \cite{BDP}
Corollaire 3 Formule (2.61), comme cons\'equence d'un r\'esultat plus
g\'en\'eral qui \'etait un d\'eveloppement complet en s\'erie
altern\'ee de la densit\'e de $S_1$ - voir \cite{BDP} Formule (2.8) :
\begin{equation}
\label{berne}f_{S_1}(x)\; =\; \sum_{n\ge 1}\frac{1}{\Gamma(\a n -1)\Gamma(1+1/\a -
  n)}x^{\a n -2}.
\end{equation} 
Quitte \`a changer de variable, ce r\'esultat est en effet \'equivalent \`a une d\'ecomposition du m\^eme type pour la densit\'e de $T_1$ :
\begin{equation}
\label{bernique}
f_{T_1}(t)\; =\; \sum_{n\ge 1}\frac{1}{\a\Gamma(\a n -1)\Gamma(1+1/\a -
  n)}t^{1/\a - n -1},
\end{equation}  
et entra\^{\i}ne donc automatiquement la Proposition 6. Pour obtenir (\ref{berne}), Bernyk, Dalang et Peskir avaient invers\'e la double
transform\'ee de Laplace (\ref{WH}) d'abord en temps et en termes de fonctions
d'erreur g\'en\'eralis\'ees, puis en espace \`a l'aide d'une \'equation
int\'egrale fractionnaire de type Abel, dont ils avaient pu donner la
solution par un d\'eveloppement en s\'erie. Remarquons au passage que les
fonctions de Mittag-Leffler d'indice plus grand que 1 interviennent dans la
r\'esolution d'\'equations int\'egrales fractionnaires de type Abel \cite{GM}. Mais il ne semble pas que l'on puisse d\'eduire de \cite{GM, BDP} l'identit\'e en loi du Th\'eor\`eme 3.

La m\'ethode de cette note, o\`u l'on a invers\'e (\ref{WH}) d'abord en espace
puis en temps, est plus simple que celle de \cite{BDP}. En revanche, l'expression de la densit\'e de $T_1$ obtenue par le Th\'eor\`eme 3 comme une
int\'egrale double :
\begin{eqnarray*}
f_{T_1}(t)& = & \frac{-\sin
  \pi\a}{\a\pi}\int_0^\infty f_{\hT_1}(t/u) \frac{1 + u^{1/\a}}{u^3 - 2u^2\cos
  \pi\a + u} du\\
& = & \frac{(-\sin \pi\a)}{\pi^2\a}\int_0^\infty \lpa \int_0^\infty e^{-(tv/u
  + v^{1/\a}\cos(\pi/\a))}\sin(v^{1/\a}\sin(\pi/\a))\, dv\rpa
\frac{(1+u^{1/\a})}{u^3 - 2u^2\cos \pi\a +u}du
\end{eqnarray*}
ne permet pas de retrouver la d\'ecomposition (\ref{bernique}) puisque le
terme $-v^{1/\a}\cos(\pi/\a) > 0$ sous l'exponentielle interdit d'appliquer
Fubini. D'un autre c\^ot\'e, on peut trouver cette derni\`ere expression un
tantinet plus ramass\'ee que celle obtenue en changeant la variable dans la formule
(2.56) de \cite{BDP}. 

Une autre formulation de la densit\'e de $T_1$ peut \^etre obtenue en inversant
(\ref{L2}) comme dans le Th\'eor\`eme 1. Avec les m\^emes notations on peut en effet \'ecrire
$$\EE[e^{-x T_1}]\; =\; \Ea(x) - \a x^{1-1/\a} \Ea'(x)\; =\; \frac{1}{2\pi
  \ii}\int_H e^t \lpa\frac{t^{\a -1} - x^{1 -1/\a}}{t^\a - x}\rpa dt.$$
D'autre part, il est possible de choisir $H = H_\delta$ telle que son
intersection avec ${\rm Re} (t) \ge 0$ soit contenue dans le cercle centr\'e
en l'origine de rayon $x^{1/\a},$ et son intersection avec ${\rm Re} (t) \le 0$ dans un secteur angulaire $-\pi/2 - \delta
<\theta < \pi/2 +\theta$ avec $\delta$ aussi petit qu'on veut - voir \`a
nouveau \cite{D} Figure 56. Comme $\a \in ]1,2[$  le changement de variable
$t\mapsto t^\a$ transforme alors $H_\delta$ en un chemin de Hankel situ\'e
dans le demi-plan ${\rm Re} (t) < x$ pour $\delta$ assez petit. Raisonnant exactement comme
dans \cite{P2} on peut donc calculer
\begin{eqnarray*}
\frac{1}{2\pi \ii}\int_H \frac{t^{\a -1} e^t}{t^\a - x} dt & = &
\frac{1}{2\pi \a\ii}\int_H \frac{e^{t^{1/\a}}}{t - x} dt \\
& = & \frac{1}{\a}\int_0^\infty e^{-xt} g_{\hT_1} (t) dt
\end{eqnarray*}
o\`u l'on a not\'e
\begin{equation}
\label{gege}
g_{\hT_1}(t)\; =\;  \frac{-1}{\pi}\int_0^\infty
e^{-(tu-u^{1/\a}\cos(\pi/\a))}\sin(u^{1/\a}\sin(\pi/\a))\, du.
\end{equation}
Le calcul du deuxi\`eme terme se fait comme plus haut et nous laissons les
d\'etails au lecteur :
\begin{eqnarray*}
\frac{1}{2\pi \ii}\int_H \frac{x^{1-1/\a} e^t}{t^\a - x} dt & = & \frac{-\sin
  \pi\a}{\pi}\int_0^\infty e^{-x^{1/\a}u} \lpa \frac{u^\a}{u^{2\a} - 2u^\a\cos
  \pi\a  + 1}\rpa du\\
& = & \int_0^\infty e^{-xt} \lpa \frac{-\sin
  \pi\a}{\a\pi}\int_0^\infty f_{\hT_1}(t/u) \frac{u^{1/\a}}{u^3 - 2u^2\cos
  \pi\a  + u} du \rpa dt,
\end{eqnarray*}
d'o\`u l'on d\'eduit 
\begin{equation}
\label{convex}
f_{T_1}(t)\; =\; (1/\a) g_{\hT_1}(t) \; + \;  (1- 1/\a) h_{\hT_1}(t)
\end{equation}
en ayant not\'e
$$h_{\hT_1}(t)\; =\;\frac{-\sin
  \pi\a}{(\a-1)\pi}\int_0^\infty f_{\hT_1}(t/u) \lpa\frac{u^{1/\a}}{u^2 - 2u\cos
  \pi\a  + 1}\rpa \frac{du}{u}\cdot$$
Puisque
$$u\;\mapsto\;\frac{(-\sin
  \pi\a)u^{1/\a}}{(\a-1)\pi(u^2 - 2u\cos
  \pi\a  + 1)}$$
est une densit\'e de probabilit\'e sur $\rl^+,$ si l'on note ${\bar T}$ la
v.a. correspondante on voit que $h_{\hT_1}$ est la densit\'e du produit
ind\'ependant ${\bar T} \times \hT_1.$ D'autre part, en comparant
(\ref{convex}) avec l'expression de $f_{T_1}$ donn\'ee par le Th\'eor\`eme 3,
on voit que $g_{\hT_1}$ est aussi une densit\'e de probabilit\'e, celle du produit
ind\'ependant ${\tilde T} \times \hT_1$ o\`u ${\tilde T}$ a pour densit\'e
$$u\;\mapsto\;\frac{-\sin
  \pi\a}{\pi(u^2 - 2u\cos
  \pi\a  + 1)}\cdot$$
En particulier $g_{\hT_1}$ est positive, ce qui
n'\'etait pas \'evident a priori, et la d\'ecomposition (\ref{convex}) se lit
comme une {\em combinaison convexe} de densit\'es de probabilit\'e. Le premier terme $(1/\a) g_{\hT_1}$ donne le comportement de $f_{T_1}$ en 0, et le second $(1- 1/\a) h_{\hT_1}$ celui de $f_{T_1}$ \`a l'infini.

Remarquons enfin que la densit\'e $g_{\hT_1}$ a quelque parent\'e avec $f_{\hT_1}$ - d'o\`u sa souscription en
$\hT_1$ - puisque son d\'eveloppement en s\'erie altern\'ee, que l'on tire de (\ref{gege}) exactement comme
dans \cite{P1}, est
$$g_{\hT_1}(t)\; =\; \sum_{n\ge
  1}\frac{1}{\Gamma (-n/\a) n!}t^{-1 -n/\a},$$
tandis que celui de $f_{\hT_1}$ est
$$f_{\hT_1}(t)\; =\; \sum_{n\ge
  1}\frac{(-1)^n}{\Gamma (-n/\a) n!}t^{-1 -n/\a}.$$
En changeant une derni\`ere fois la variable, on trouve donc l'expression d'une
certaine transform\'ee int\'egrale reliant ces deux s\'eries altern\'ees :
$$ g_{\hT_1}(t)\; =\;\int_0^\infty f_{\hT_1}(u) \lpa\frac{(-\sin
  \pi\a)u}{\pi(u^2 - 2\cos
  \pi\a tu + t^2)}\rpa du$$
avec laquelle je concluerai, en remerciant Jean-Fran\c{c}ois Burnol et Francis Hirsch pour leur aide dans l'\'elaboration de cet article.

\end{document}